\documentclass[11pt,reqno]{amsproc}

\title[]{Remarks on the fractional Laplacian with Dirichlet boundary conditions and applications}

\author{Peter Constantin}
\address{Department of Mathematics, Princeton University, Princeton, NJ 08544}
\email{const@math.princeton.edu}

\author{Mihaela Ignatova}
\address{Department of Mathematics, Princeton University, Princeton, NJ 08544}
\email{ignatova@math.princeton.edu}

\usepackage[margin=1in]{geometry}
\usepackage{amsmath, amsthm, amssymb}
\usepackage{times}
\usepackage{color}
\usepackage{hyperref}

\newcommand{\pa}{\partial}
\newcommand{\la}{\label}
\newcommand{\fr}{\frac}
\newcommand{\na}{\nabla}
\newcommand{\be}{\begin{equation}}
\newcommand{\ee}{\end{equation}}
\newcommand{\ba}{\begin{array}{l}}
\newcommand{\ea}{\end{array}}
\newcommand{\Rr}{{\mathbb R}}

\newcommand{\beg}{\begin}

\newcommand{\G}{\Gamma}

\newcommand{\D}{\Delta}
\renewcommand{\l}{\Lambda_D}
\date{today}
\begin{document}
\begin{abstract}
We prove nonlinear lower bounds and commutator estimates for the Dirichlet fractional Laplacian in bounded domains. The applications include bounds for linear drift-diffusion equations with nonlocal dissipation and global existence of weak solutions of critical surface quasi-geostrophic equations. 
\end{abstract}

\maketitle
\section{Introduction}
Drift-diffusion equations with nonlocal dissipation naturally occur in hydrodynamics and in models of electroconvection. The study of these equations in bounded domains is hindered by a lack of explicit information on the kernels of the nonlocal operators appearing in them. In this paper we develop tools adapted for the Dirichlet boundary case:  the C\'{o}rdoba-C\'{o}rdoba inequality (\cite{cc}) and a nonlinear lower bound in the spirit of (\cite{cv1}), and commutator estimates. Lower bounds for the fractional Laplacian are instrumental in proofs of regularity of solutions to nonlinear nonlocal drift-diffusion equations. The presence of boundaries requires natural modifications of the bounds. The nonlinear bounds are proved using a representation based on the heat kernel and fine information regarding it (\cite{davies1}, \cite{qszhang1}, \cite{qszhang2}). Nonlocal diffusion operators in bounded domains do not commute in general with differentiation. The commutator estimates are proved  using the method of harmonic extension and results of (\cite{cabre}). We apply these tools to linear drift-diffusion equations with nonlocal dissipation, where we obtain strong global bounds, and to global existence of weak solutions of the surface quasi-geostrophic equation (SQG) in bounded domains.  

We consider a bounded open domain $\Omega\subset \Rr^d$ with smooth (at least $C^{2,\alpha}$) boundary. We denote by $\D$ the Laplacian operator with homogeneous Dirichlet boundary conditions. Its $L^2(\Omega)$ - normalized eigenfunctions are denoted $w_j$, and its eigenvalues counted with their multiplicities are denoted $\lambda_j$: 
\be
-\D w_j = \lambda_j w_j.
\la{ef}
\ee
It is well known that $0<\lambda_1\le...\le \lambda_j\to \infty$  and that $-\D$ is a positive selfadjoint operator in $L^2(\Omega)$ with domain ${\mathcal{D}}\left(-\D\right) = H^2(\Omega)\cap H_0^1(\Omega)$.
The ground state $w_1$ is positive and
\be
c_0d(x) \le w_1(x)\le C_0d(x)
\la{phione}
\ee
holds for all $x\in\Omega$, where 
\be
d(x) = dist(x,\pa\Omega)
\la{dx}
\ee
and $c_0, \, C_0$ are positive constants depending on $\Omega$. Functional calculus can be defined using the eigenfunction expansion. In particular
\be
\left(-\D\right)^{\alpha}f = \sum_{j=1}^{\infty}\lambda_j^{\alpha} f_j w_j
\la{funct}
\ee
with 
\[
f_j =\int_{\Omega}f(y)w_j(y)dy
\]
for $f\in{\mathcal{D}}\left(\left (-\D\right)^{\alpha}\right) = \{f\left |\right. \; (\lambda_j^{\alpha}f_j)\in \ell^2(\mathbb N)\}$.
We will denote by
\be
\l^s = \left(-\D\right)^{\alpha}, \quad s=2\alpha
\la{lambdas}
\ee
the fractional powers of the Dirichlet Laplacian, with $0\le \alpha\le 1$
and with $\|f\|_{s,D}$ the norm in ${\mathcal{D}}\left (\l^s\right)$:
\be
\|f\|_{s,D}^2 = \sum_{j=1}^{\infty}\lambda_j^{s}f_j^2.
\la{norms}
\ee
It is well-known and easy to show that
\[
{\mathcal{D}}\left( \l \right) = H_0^1(\Omega).
\]
Indeed, for $f\in{\mathcal{D}}\left (-\D\right)$ we have
\[
\|\na f\|^2_{L^2(\Omega)} = \int_{\Omega}f\left(-\D\right)fdx = \|\l f\|_{L^2(\Omega)}^2 = \|f\|^2_{1,D}. 
\]
We recall that the Poincar\'{e} inequality implies that the Dirichlet integral on the left-hand side above is equivalent to the norm in $H_0^1(\Omega)$
and therefore the identity map from the dense subset ${\mathcal{D}}\left(-\D\right)$ of $H_0^1(\Omega)$ to ${\mathcal D}\left(\l\right)$ is an isometry, and thus $H_0^1(\Omega)\subset {\mathcal{D}}\left(\l\right)$. But ${\mathcal{D}}\left(-\D\right)$ is dense in  ${\mathcal D}\left(\l\right)$ as well, because finite linear combinations of eigenfunctions are dense in  ${\mathcal D}\left(\l\right)$. Thus the opposite inclusion is also true, by the same isometry argument.  \\
Note that in view of the identity
\be
\lambda^{\alpha} = c_{\alpha}\int_0^{\infty}(1-e^{-t\lambda})t^{-1-\alpha}dt,
\la{lambdalpha}
\ee
with 
\[
1 = c_{\alpha} \int_0^{\infty}(1-e^{-s})s^{-1-\alpha}ds,
\]
valid for $0\le \alpha <1$, we have the representation
\be
\left(\left(-\D\right)^{\alpha}f\right)(x) = c_{\alpha}\int_0^{\infty}\left[f(x)-e^{t\D}f(x)\right]t^{-1-\alpha}dt
\la{rep}
\ee
for $f\in{\mathcal{D}}\left(\left (-\D\right)^{\alpha}\right)$.
We use precise upper and lower bounds for the kernel $H_D(t,x,y)$ of the heat operator,
\be
(e^{t\D}f)(x) = \int_{\Omega}H_D(t,x,y)f(y)dy .
\la{heat}
\ee
These are as follows (\cite{davies1},\cite{qszhang1},\cite{qszhang2}).
There exists a time $T>0$ depending on the domain $\Omega$ and constants
$c$, $C$, $k$, $K$, depending on $T$ and $\Omega $ such that
\be
\ba
c\min\left (\fr{w_1(x)}{|x-y|}, 1\right)\min\left (\fr{w_1(y)}{|x-y|}, 1\right)t^{-\fr{d}{2}}e^{-\fr{|x-y|^2}{kt}}\le \\H_D(t,x,y)\le C
\min\left (\fr{w_1(x)}{|x-y|}, 1\right)\min\left (\fr{w_1(y)}{|x-y|}, 1\right)t^{-\fr{d}{2}}e^{-\fr{|x-y|^2}{Kt}}
\ea
\la{hb}
\ee
holds for all $0\le t\le T$. Moreover
\be
\fr{\left |\na_x H_D(t,x,y)\right|}{H_D(t,x,y)}\le
C\left\{
\ba
\fr{1}{d(x)},\quad\quad \quad\quad {\mbox{if}}\; \sqrt{t}\ge d(x),\\
\fr{1}{\sqrt{t}}\left (1 + \fr{|x-y|}{\sqrt{t}}\right),\;{\mbox{if}}\; \sqrt{t}\le d(x)
\ea
\right.
\la{grbx}
\ee
holds for all $0\le t\le T$.
Note that, in view of
\be
H_D(t,x,y) = \sum_{j=1}^{\infty}e^{-t\lambda_j}w_j(x)w_j(y) ,
\la{hphi}
\ee
elliptic regularity estimates and Sobolev embedding which imply uniform absolute convergence of the series (if $\pa\Omega$ is smooth enough), we have that
\be
\pa_1^{\beta}H_D(t,y,x) = \pa_2^{\beta}H_D(t,x,y)
= \sum_{j=1}^{\infty}e^{-t\lambda_j}\pa_y^{\beta}w_j(y)w_j(x)
\la{dh}
\ee
for positive $t$, where we denoted by $\pa_1^{\beta}$ and $\pa_2^{\beta}$ derivatives with respect to the first spatial variables and the second spatial variables, respectively.

Therefore, the gradient bounds (\ref{grbx}) result in
\be
\fr{\left |\na_y H_D(t,x,y)\right|}{H_D(t,x,y)}\le
C\left\{
\ba
\fr{1}{d(y)},\quad\quad \quad\quad\quad {\mbox{if}}\; \sqrt{t}\ge d(y),\\
\fr{1}{\sqrt{t}}\left (1 + \fr{|x-y|}{\sqrt{t}}\right),\;{\mbox{if}}\; \sqrt{t}\le d(y).
\ea
\right.  
\la{grby}
\ee

\section{The C\'{o}rdoba - C\'{o}rdoba inequality}
\beg{prop}{\la{cordoba}}
Let $\Phi$ be a $C^2$ convex function satisfying $\Phi(0)= 0$. Let $f\in C_0^{\infty}(\Omega)$ and let $0\le s\le 2$. Then
\be
\Phi'(f)\l^s f - \l^s(\Phi(f))\ge 0
\la{cor}
 \ee
holds pointwise almost everywhere in $\Omega$.
\end{prop}
\noindent{\bf Proof.} 
In view of the fact that both $f\in H_0^1(\Omega)\cap H^2(\Omega)$ and
$\Phi(f)\in H_0^1(\Omega)\cap H^2(\Omega)$, the terms in the inequality (\ref{cor}) are well defined. We
define
\be
\left[\left(-\D\right)^{\alpha}f\right]_{\epsilon}(x)= c_{\alpha}\int_\epsilon^{\infty}\left[f(x)-e^{t\D}f(x)\right]t^{-1-\alpha}dt
\ee
and approximate the representation (\ref{rep}):
\be
\left(\left(-\D\right)^{\alpha}f\right)(x) = \lim_{\epsilon\to 0}\left[\left(-\D\right)^{\alpha}f\right]_{\epsilon}(x).
\la{epsap}
\ee
The limit is strong in $L^2(\Omega)$. We start the calculation with this approximation and then we rearrange terms:
\[
\ba
\Phi'(f(x))\left[\l^{2\alpha} f\right]_{\epsilon}(x) - \left[\l^{2\alpha} (\Phi(f))\right]_{\epsilon}(x) \\ 
= c_{\alpha}\int_\epsilon^{\infty} t^{-1-\alpha}dt\int_{\Omega}\left\{\Phi'(f(x))\left[\fr{1}{|\Omega|}f(x) -H_D(t,x,y)f(y)\right]
-\fr{1}{|\Omega|}\Phi(f(x)) + H_D(t,x,y)\Phi(f(y))\right\}dy\\
=c_{\alpha}\int_\epsilon^{\infty} t^{-1-\alpha}dt\int_{\Omega}H_D(t,x,y)\left[\Phi(f(y))-\Phi(f(x)) - \Phi'(f(x))(f(y)-f(x))\right]dy\\
+ c_{\alpha}\int_\epsilon^{\infty}t^{-1-\alpha}dt\int_{\Omega}\left[f(x)\Phi'(f(x))-\Phi(f(x))\right](\fr{1}{|\Omega|} - H_D(t,x,y))dy\\
= c_{\alpha}\int_\epsilon^{\infty} t^{-1-\alpha}dt\int_{\Omega}H_D(t,x,y)\left[\Phi(f(y))-\Phi(f(x)) - \Phi'(f(x))(f(y)-f(x))\right]dy\\
+\left[f(x)\Phi'(f(x))-\Phi(f(x))\right]c_{\alpha}\int_\epsilon^{\infty}t^{-1-\alpha}(1-e^{t\D}1)dt 
\ea
\]
Because of the convexity of $\Phi$ we have
\[
\Phi(b)-\Phi(a) -\Phi'(a)(b-a)\ge 0, \quad \forall\;\;a, b\in \Rr,
\]
and because $\Phi(0) =0$ we have 
\[
a\Phi'(a)\ge \Phi(a), \quad \forall \;\; a\in\Rr.
\]
Consequently $f(x)\Phi'(f(x))-\Phi(f(x))\ge 0$ holds everywhere. The function
\[
\theta = e^{t\D}1
\]
solves the heat equation $\pa_t\theta -\D\theta =0$ in $\Omega$, with homogeneous Dirichlet boundary conditions, and with initial data equal everywhere to $1$. Although $1$ is not in the domain of $-\D$, $e^{t\D}$ has a unique extension to $L^2(\Omega)$ where $1$ does belong, and on the other hand, by the maximum principle $0\le\theta(x,t)\le 1$ holds for $t\ge 0$, $x\in\Omega$. It is only because $1\notin{\mathcal{D}}(-\D)$ that we had to use the $\epsilon$ approximation. Now we discard the nonnegative term 
\[
\left[f(x)\Phi'(f(x))-\Phi(f(x))\right]c_{\alpha}\int_{\epsilon}^{\infty}(1-\theta(x,t))t^{-1-\alpha}dt
\] in the calculation above, and deduce that
\be
\Phi'(f(x))\left[\l^{2\alpha} f\right]_{\epsilon}(x) - \left[\l^{2\alpha} (\Phi(f))\right]_{\epsilon}(x) \ge 0
\ee
as an element of $L^2(\Omega)$. (This simply means that its integral against any nonnegative $L^2(\Omega)$ function is nonnegative.) Passing to the limit $\epsilon\to 0$ we obtain the inequality (\ref{cor}). If $\Phi $ and the boundary of the domain are smooth enough then we can prove that the terms in the inequality are continuous, and therefore the inequality holds everywhere.

\section{The Nonlinear Bound}
We prove a bound in the spirit of (\cite{cv1}). The nonlinear lower bound was
used as an essential ingredient in proofs of global regularity for drift-diffusion equations with nonlocal dissipation.
\beg{thm}
Let $f\in L^{\infty}(\Omega)\cap {\mathcal D}(\l^{2\alpha})$, $0\le \alpha<1$. 
Assume that $f= \pa q$ with $q \in L^{\infty}(\Omega)$ and $\pa$ a first order derivative. Then there exist constants $c$, $C$ depending on $\Omega$ and $\alpha$ such that
\be
f\l^{2\alpha} f -\fr{1}{2}\l^{2\alpha} f^2 \ge c\|q\|_{L^{\infty}}^{-2\alpha} |f_d|^{2+2\alpha}
\la{nlb}
\ee
holds pointwise in $\Omega$, with
\be
|f_d(x)| = \left\{
\ba
|f(x)|,\quad\quad {\mbox{if}}\;\; |f(x)| \ge C\|q\|_{L^{\infty}(\Omega)}\max\left(\fr{1}{\rm diam(\Omega)},\fr{1}{d(x)}\right),\\
0,\quad\quad\quad\quad {\mbox{if}}\;\; |f(x)| \le C\|q\|_{L^{\infty}(\Omega)}\max\left(\fr{1}{\rm diam(\Omega)},\fr{1}{d(x)}\right).
\ea
\right.
\la{fd}
\ee
\end{thm}
\noindent{\bf{Proof.}} We start the calculation using the inequality 
\be
f\l^{2\alpha}f - \fr{1}{2}\l^{2\alpha} f^2 \ge \fr{1}{2}c_{\alpha}\int_0^{\infty}\psi\left(\fr{t}{\tau}\right)t^{-1-\alpha}dt\int_{\Omega}H_D(t,x,y)(f(x)-f(y))^2dy
\la{lowone}
\ee
where $\tau>0$ is arbitrary and $0\le \psi(s)\le 1$ is a smooth function, vanishing identically for $0\le s\le 1$ and equal identically to $1$ for $s\ge 2$. This follows repeating the calculation of the proof of the C\'{o}rdoba-C\'{o}rdoba inequality with $\Phi(f)=\fr{1}{2}f^2$:
\[
\ba
f(x)\left[\l^{2\alpha}f\right]_{\epsilon}(x) - \fr{1}{2}\left[\l^{2\alpha} f^2\right]_{\epsilon}(x) \\
= c_{\alpha}\int_\epsilon^{\infty}t^{-1-\alpha}\int_{\Omega}\left\{\left[\fr{1}{|\Omega|} f(x)^2 - f(x)H_D(t,x,y)f(y)\right]- \fr{1}{2|\Omega|}f^{2}(x) + \fr{1}{2}H_D(t,x,y)f^2(y)\right\}dy\\
=c_{\alpha}\int_\epsilon^{\infty}t^{-1-\alpha}dt\int_{\Omega}\left\{\fr{1}{2}\left[H_D(t,x,y)(f(x)-f(y))^2\right]  + \fr{1}{2}f^2(x)\left[\fr{1}{|\Omega|} -H_D(t,x,y)\right]\right\}dy \\
= c_{\alpha}\int_\epsilon^{\infty}t^{-1-\alpha}dt\int_{\Omega}\left\{\fr{1}{2}\left[H_D(t,x,y)(f(x)-f(y))^2\right]dy + \fr{1}{2}f^2(x)\left[1-e^{t\D}1\right](x)\right\}\\
\ge c_{\alpha}\int_\epsilon^{\infty}t^{-1-\alpha}dt\int_{\Omega}\fr{1}{2}H_D(t,x,y)\left(f(x)-f(y)\right)^2dy  
\ea
\]
where in the last inequality we used the maximum principle again. Then, we choose $\tau>0$ and let $\epsilon<\tau$. It follows that
\[
f(x)\left[\l^{2\alpha}f\right]_{\epsilon}(x) - \fr{1}{2}\left[\l^{2\alpha} f^2\right]_{\epsilon}(x) \ge
\fr{1}{2}c_{\alpha}\int_0^{\infty}\psi\left(\fr{t}{\tau}\right)t^{-1-\alpha}dt\int_{\Omega}H_D(t,x,y)\left(f(x)-f(y)\right)^2dy.
\]
We obtain (\ref{lowone}) by letting $\epsilon\to 0$. We restrict to $t\le T$,
\be
\left[f\l^{2\alpha}f - \fr{1}{2}\l^{2\alpha} f^2\right](x) \ge
\fr{1}{2}c_{\alpha}\int_0^{T}\psi\left(\fr{t}{\tau}\right)t^{-1-\alpha}dt\int_{\Omega}H_D(t,x,y)\left(f(x)-f(y)\right)^2dy
\la{lowtwo}
\ee
and open brackets in (\ref{lowtwo}):
\be
\ba
\left[f\l^{2\alpha}f - \fr{1}{2}\l^{2\alpha} f^2\right](x) \\
\ge \fr{1}{2}f^2(x)c_{\alpha}\int_0^T\psi\left(\fr{t}{\tau}\right)t^{-1-\alpha}dt\int_{\Omega}H_D(t,x,y)dy 
- f(x)c_{\alpha}\int_0^T\psi\left(\fr{t}{\tau}\right)t^{-1-\alpha}dt\int_{\Omega}H_D(t,x,y)f(y)dy\\
\ge |f(x)|\left [ \fr{1}{2}|f(x)| I(x) - J(x)\right]
\ea
\la{lowthree}
\ee
with
\be
I(x) = c_{\alpha}\int_0^T\psi\left(\fr{t}{\tau}\right)t^{-1-\alpha}dt\int_{\Omega}H_D(t,x,y)dy,
\la{ix}
\ee
and
\be
\ba
J(x) = c_{\alpha}\left |\int_0^T\psi\left(\fr{t}{\tau}\right)t^{-1-\alpha}dt\int_{\Omega}H_D(t,x,y)f(y)dy\right |\\
= c_{\alpha}\left |\int_0^T\psi\left(\fr{t}{\tau}\right)t^{-1-\alpha}dt\int_{\Omega}\pa_yH_D(t,x,y)q(y)dy\right |.
\ea
\ee
We proceed with a lower bound on $I$ and an upper bound on $J$.  For the lower bound on $I$ we note that
\[
\theta (x,t) = \int_{\Omega}H_D(t,x,y)dy\ge \int_{|x-y|\le \fr{d(x)}{2}} H_D(t,x,y)dy
\]
because $H_D$ is positive. Using the lower bound in (\ref{phione}) we have that
$|x-y|\le \fr{d(x)}{2}$ implies
\[
\fr{w_1(x)}{|x-y|}\ge 2c_0,\quad \fr{w_1(y)}{|x-y|}\ge c_0,
\]
and then, using the lower bound in (\ref{hb}) we obtain
\[
H_D(t,x,y) \ge 2cc_0^2t^{-\fr{d}{2}}e^{-\fr{|x-y|^2}{kt}}. 
\]
Integrating it follows that
\[
\theta(x,t) \ge 2 cc_0^2\omega_{d-1}k^{\fr{d}{2}}\int_0^{\fr{d(x)}{2\sqrt{kt}}}\rho^{d-1}e^{-\rho^2}d\rho
\]
If $\fr{d(x)}{2\sqrt{kt}}\ge 1$ then the integral is bounded below by
$\int_0^1\rho^{d-1}e^{-\rho^2}d\rho$. If $\fr{d(x)}{2\sqrt{kt}}\le 1$ then
$\rho\le 1$ implies that the exponential is bounded below by $e^{-1}$ and so
\be
\theta(x,t)\ge c_1\min\left\{1, \left(\fr{d(x)}{\sqrt{t}}\right)^d\right\}
\la{thetalow}
\ee
for all $0\le t\le T$ where $c_1$ is a positive constant, depending on $\Omega$. Because 
\[
I(x) = \int_0^T\psi\left(\fr{t}{\tau}\right)t^{-1-\alpha}\theta(x,t)dt
\]
we have
\[
\ba
I(x)\ge c_1\int_0^{\min(T, d^2(x))}\psi\left(\fr{t}{\tau}\right)t^{-1-\alpha}dt\\
= c_1\tau^{-\alpha}\int_1^{\tau^{-1}(\min(T, d^2(x)))}\psi(s)s^{-1-\alpha}ds
\ea
\]
Therefore we have that
\be
I(x)\ge c_2 \tau^{-\alpha}
\la{ilow}
\ee
with $c_2 = c_1\int_1^2\psi(s)s^{-1-\alpha}ds$, a positive constant depending only on $\Omega$ and $\alpha$, provided $\tau$ is small enough,
\be
\tau \le \fr{1}{2}\min(T,  d^2(x)).
\la{taucond}
\ee
In order to bound $J$ from above we use the upper bound (\ref{grby}) which yields
\be
\int_{\Omega}|\na_y H_D(t,x,y)|dy \le C_1 t^{-\fr{1}{2}}
\la{grup}
\ee
with $C_1$ depending only on $\Omega$. Indeed, 
\[
\ba
\int_{d(y)\ge \sqrt{t}}|\na_y H_D(t,x,y)|dy \\\le C_2t^{-\fr{1}{2}}
\int_{\Rr^d}\left (1 + \fr{|x-y|}{\sqrt{t}}\right)t^{-\fr{d}{2}}e^{-\fr{|x-y|^2}{kt}}dy\\
= C_3t^{-\fr{1}{2}}
\ea
\]
and, in view of the upper bound in (\ref{phione}), $\fr{1}{d(y)}w_1(y)\le C_0$ and the upper bound in (\ref{hb}),
\[
\ba
\int_{d(y)\le \sqrt{t}}|\na_y H_D(t,x,y)|dy \\\le C_4\int_{\Rr^d}\fr{1}{|x-y|}t^{-\fr{d}{2}}e^{-\fr{|x-y|^2}{Kt}}dy = C_5t^{-\fr{1}{2}}
\ea
\]
Now 
\[
J\le \|q\|_{L^{\infty}(\Omega)}\int_0^T\psi\left(\fr{t}{\tau}\right)t^{-1-\alpha}dt\int_{\Omega}|\na_y H_D(t,x,y)|dy
\]
and therefore, in view of (\ref{grup})
\[
J\le C_1 \|q\|_{L^{\infty}(\Omega)}\int_0^T\psi\left(\fr{t}{\tau}\right )t^{-\fr{3}{2}-\alpha}dt
\]
and therefore
\be
J \le C_6\|q\|_{L^{\infty}(\Omega)}\tau^{-\fr{1}{2}-\alpha}
\la{jup}
\ee
with
\[
C_6 = C_1\int_1^{\infty}\psi(s)s^{-\fr{3}{2}-\alpha}ds
\]
a constant depending only on $\Omega$ and $\alpha$. Now, because of the lower bound (\ref{lowthree}), 
if we can choose $\tau$ so that
\[
J(x) \le \fr{1}{4} |f(x)|I(x)
\]
then it follows that
\be
\left[f\l^{2\alpha}f - \fr{1}{2}\l^{2\alpha} f^2\right](x) \ge \fr{1}{4}f^2(x)I(x).
\la{lowfour}
\ee
Because of the bounds (\ref{ilow}), (\ref{jup}) the choice
\be
\tau(x) = c_3\fr{\|q\|_{L^{\infty}}^2}{|f(x)|^2}
\la{tauchoice}
\ee
with $c_3= 16 C_6^2c_2^{-2}$ achieves the desired bound. The requirement (\ref{taucond}) limits the possibility of making this choice to the situation
\be
c_3\fr{\|q\|_{L^{\infty}}^2}{|f(x)|^2} \le \fr{1}{2}\min(T,  d^2(x))
\ee
which leads to the statement of the theorem. Indeed, if (\ref{tauchoice}) is allowed then the lower bound in (\ref{lowfour}) becomes
\be
\left[f\l^{2\alpha}f - \fr{1}{2}\l^{2\alpha} f^2\right](x) \ge  c\|q\|_{L^{\infty}}^{-2\alpha} |f_d|^{2+2\alpha}
\la{lowfive}
\ee
with $c=\fr{1}{4}c_2c_3^{-\alpha}$. 
\section{Commutator estimates}
We start by considering the commutator $[\na,\l]$ in $\Omega = \Rr^d_+$.  The heat kernel with Dirichlet boundary conditions is
\[
H(x,y,t) = ct^{-\fr{d}{2}}\left( e^{-\fr{|x-y|^2}{4t}} - e^{-\fr{|x-{\widetilde {y}}|^2}{4t}}\right) 
\]
where $\widetilde{y} = (y_1,\dots, y_{d-1}, -y_d)$. We claim that
\be
\int_{\Omega} (\na_x+ \na_y)H(x,y,t) dy \le Ct^{-\fr{1}{2}}e^{-\fr{x_d^{2}}{4t}}.
\la{hone}
\ee
Indeed, the only nonzero component occurs when the differentiation is with respect to the normal direction, and then
\[
(\pa_{x_d} + \pa_{y_d})H(x,y,t) = ct^{-\fr{d}{2}}e^{-\fr{|x'-y'|^2}{4t}}\left(\fr{x_d+y_d}{t}\right) e^{-\fr{(x_d+y_d)^2}{4t}}
\] 
where we denoted $x' = (x_1,\dots, x_{d-1})$ and $y'=(y_1, \dots, y_{d-1})$.
Therefore
\[
\ba
\int_{\Omega} (\na_x+ \na_y)H(x,y,t) dy \le Ct^{-\fr{1}{2}}\int_0^{\infty}\left(\fr{x_d+y_d}{t}\right )e^{-\fr{(x_d+y_d)^2}{4t}}dy_d\\
= Ct^{-\fr{1}{2}}\int_{\fr{x_d}{\sqrt {t}}}^{\infty}\xi e^{-\fr{\xi^2}{4}}d\xi
\\ 
= Ct^{-\fr{1}{2}}e^{-\fr{x_d^2}{4t}}.
\ea
\]
Consequently
\[
K(x,y) = \int_0^{\infty}t^{-\fr{3}{2}}(\na_x + \na_y) H(x,y,t)dt
\]
obeys
\[
\int_{\Omega}K(x,y) dy \le C\int_0^{\infty}t^{-2}e^{-\fr{x_d^2}{4t}}dt 
= \fr{C}{x_{d}^{2}}.
\]
The commutator $[\na, \l]$ is computed as follows 
\[
\ba
[\na, \l]f(x) = \int_0^{\infty}t^{-\fr{3}{2}}\int_{\Omega}\left[\na_xH_D(x,y,t)f(y) - H_D(x,y,t)\na_y f(y)\right]dydt\\
=\int_0^{\infty}t^{-\fr{3}{2}}\int_{\Omega}(\na_x+\na_y)H_{D}(x,y,t)f(y)dydt \\
=\int_{\Omega}K(x,y)f(y)dy.
\ea
\]
We have proved thus that the kernel $K(x,y)$ of the commutator obeys
\be
\int_{\Omega}K(x,y) dy \le Cd(x)^{-2}
\la{kernelone}
\ee
and therefore we obtain, for instance, for any $p,q\in [1,\infty]$ with $p^{-1}+q^{-1}=1$
\[
\left |\int_{\Omega} g[\na,\l]fdx\right | \le C\left(\int_{\Omega}d(x)^{-2}|f(x)|^pdx\right)^{\fr{1}{p}}\left(\int_{\Omega}d(x)^{-2}|g(x)|^qdx\right)^{\fr{1}{q}}.
\]

In general domains, the absence of explicit expressions for the heat kernel with Dirichlet boundary conditions requires a less direct approach to commutator estimates.

We take thus an open bounded domain $\Omega\subset\Rr^d$ with smooth boundary and describe the square root of the Dirichlet Laplacian using the harmonic extension. We denote 
\[
Q= \Omega\times\Rr_+ =\{(x,z)\left|\right. \; x\in\Omega, z>0\}
\]
and consider the traces of functions in $H_{0,L}^1(Q)$,
\[
H_{0,L}^1(Q) = \{v\in H^1(Q)\; \left |\right. \; v(x,z)=0,\; x\in\pa\Omega, \; z>0 \}
\]
\be
V_0(\Omega) = \{f\;\left |\right. \;\exists v\in H_{0,L}^1(Q),\; f(x)= v(x,0), x\in\Omega\}
\la{vzero}
\ee
where we slightly abused notation by referring to the images of $v$ under restriction operators as $v(x,z)$ for $x\in\pa\Omega$, and as $v(x,0)$ for $x\in \Omega$. 
We recall from (\cite{cabre}) that, on one hand,
\be
V_0(\Omega) = \{f\in H^{\fr{1}{2}}(\Omega)\;\left|\right. \; \int_{\Omega}\fr{f^2(x)}{d(x)}dx <\infty\}
\la{vzerochar1}
\ee
with norm
\[
\|f\|^2_{V_0}= \|f\|^2_{H^{\fr{1}{2}}(\Omega)} + \int_{\Omega}\fr{f^2(x)}{d(x)}dx,
\]
and on the other hand $V_0(\Omega) = {\mathcal D}(\l^{\fr{1}{2}})$, i.e.
\be
V_0(\Omega) = \{f\in L^2(\Omega)\;\left |\right.\; f = \sum_jf_jw_j,\; \sum_j\lambda_j^{\fr{1}{2}}f_j^2<\infty\}
\la{vzerochar2}
\ee
with equivalent norm
\[
\|f\|^2_{\fr{1}{2}, D} = \sum_{j=1}^{\infty}\lambda_j^{\fr{1}{2}}f_j^2 = \|\l^{\fr{1}{2}}f\|_{L^2(\Omega)}^2.
\]
The harmonic extension of $f$ will be denoted $v_f$. It is given by
\be
v_f(x,z) = \sum_{j=1}^{\infty}f_je^{-z\sqrt{\lambda_j}}w_j(x)
\la{vf}
\ee
and the operator $\l$ is then identified with
\be
\l f = -\left(\pa_z v_f\right)_{\left| \right. \;z=0}
\la{lambdavf}
\ee
Note that if $f\in V_0(\Omega)$ then $v_f\in H^1(Q)$.  Note also, that $v_f$ decays exponentially in the sense that
\be
\|v_f\|_{e^{zl}H^1(Q)} = \|e^{zl}\na v_f\|_{L^2(Q)} + \|e^{zl}v_f\|_{L^2(Q)} \le C \|f\|_{V_0}
\la{expb}
\ee
holds with $\ell = \fr{\lambda_1}{4}$.  We use a lemma in  $Q$:
\beg{lemma}{\la{z}} Let $F\in H^{-1}(Q)$ (the dual of $H_0^1(Q)$).
Then the problem
\be
\left\{
\ba
-\Delta u = F, \quad {\mbox{in}}\; Q,\\
u = 0, \quad {\mbox{on}}\; \pa Q
\ea
\right.
\la{delu}
\ee
has a unique weak solution $u\in H_0^1(Q)$. If $F\in L^2(Q)$ and if there exists $l>0$ so that 
\[
\|e^{zl}F\|_{L^2(Q)}^2=\int e^{2zl}|F(x,z)|^2dxdz <\infty
\]
then $u\in H_0^1(Q)\cap H^2(Q)$ and it satisfies
\[
\|u\|_{H^2(Q)}\le C\|e^{zl}F\|_{L^2(Q)}
\]
with $C$ a constant depending only on $\Omega$ and $l$. 
\end{lemma}
{\noindent{\bf Proof.}} We consider the domain $U= \Omega\times \Rr$ and take the odd extension of  $F$ to $U$, $F(x,-z)=-F(x,z)$. The existence of a weak solution in $H_0^1(U)$ follows by variational methods, by minimizing
\[
I(v) = \int_U\left(\fr{1}{2}|\na v|^2 + vF\right)dxdz
\]
among all odd functions $v\in H_0^1(U)$. The domain $U$ has finite width, so the Poincar\'{e} inequality
\[
\|\na v\|^2_{L^2(U)}\ge c\|v\|^2_{L^2(U)}
\]
is valid for functions in $H_0^1(U)$. This allows to show existence and uniqueness of weak solutions.
If $F\in L^2(U)$ we obtain locally uniform elliptic estimates
\[
\|u\|_{H^2(U_j)}\le C\|F\|_{L^2(V_{j})}
\]
where $U_j= \{(x,z)\left|\right.\, x\in\Omega, z\in (j-1,j+1)\}$, $V_j= \{(x,z)\left|\right.\, x\in\Omega, z\in (j-2,j+2)\}$, and $j =\pm\fr{1}{2}, \pm 1,\pm\fr{3}{2},\dots$, i.e. $j\in \fr{1}{2}{\mathbb Z}$. The constant $C$ does not depend on $j$. Because of the decay assumption on $F$, the
estimates can be summed. 
\beg{thm}
Let $a\in B(\Omega)$ where $B(\Omega) = W^{2,d}(\Omega)\cap W^{1,\infty}(\Omega)$, if $d\ge 3$, and $B(\Omega) = W^{2,p}(\Omega)$ with $p>2$, if $d=2$. There exists a constant $C$,
depending only on $\Omega$, such that
\be
\|[a,\l]f\|_{\fr{1}{2}, D}\le C\|a\|_{B(\Omega)}\|f\|_{\fr{1}{2}, D}
\la{comm}
\ee
holds for any $f\in V_0(\Omega)$, with
\[
\|a\|_{B(\Omega)} = \|a\|_{W^{2,d}(\Omega)} + \|a\|_{W^{1,\infty}(\Omega)}
\]
if $d\ge 3$ and
\[
\|a\|_{B(\Omega)} = \|a\|_{W^{2,p}(\Omega)}
\]
with $p>2$, if $d=2$.
\end{thm}

\noindent{\bf{Proof.}}
In order to compute $v_{af}$, let us note that $av_f\in H_{0,L}^1(Q)$, and
\[
\Delta(av_f) = v_f\Delta_x a + 2\na_x a \cdot \na v_f 
\]
and, because $v_f\in e^{zl}H^1(Q)$ and $a\in B(\Omega)$ we have that
\[
\|\Delta (av_f)\|_{L^2(e^{zl}dzdx)}\le C \|a\|_{B(\Omega)}\|v_f\|_{e^{zl}H^1(Q)}.
\]
Solving
\[
\left\{
\ba
\Delta u = \Delta (av_f) \quad {\mbox{in}}\; Q,\\
u = 0 \quad {\mbox{on}} \;\pa Q,
\ea
\right.
\]
we obtain $u\in H_0^1(Q)\cap H^2(Q)$. This follows from Lemma \ref{z} above. Note that $\pa_z u\in H_{0,L}^1(Q)$. Then
\[
v_{af} = av_f -u
\]
and
\[
a\l f-\l(af) = -a(\pa_z v_f)_{\left|\right. z=0} +\pa_z(av_f-u)_{\left|\right. z=0}= -\pa_z u_{\left |\right.z=0}.
\]
The estimate follows from elliptic estimates and restriction estimates
\[
\|\pa_zu_{\left|\right. z=0}\|_{V_0}\le C\|\pa_z u\|_{H^1(Q)}\le C\|a\|_{B(\Omega)}\|v_f\|_{e^{zl}H^1(Q)}\le C\|a\|_{B(\Omega)}\|f\|_{V_0}
\]

\beg{thm}{\la{comthm}} Let a vector field $a$ have components in $B(\Omega)$ defined above, $a\in \left(B(\Omega)\right)^d$. Assume that the normal component of the trace of $a$ on the boundary vanishes, 
\[
a_{\left |\right. \pa \Omega}\cdot n = 0
\]
(i.e the vector field is tangent to the boundary). 
There exists a constant $C$ such that
\be
\|[a\cdot\na ,\l]f\|_{\fr{1}{2}, D} \le C\|a\|_{B(\Omega)}\|f\|_{\fr{3}{2},D}
\la{comtwo}
\ee
holds for any $f$ such that $f\in {\mathcal{D}}\left(\l^{\fr{3}{2}}\right)$.
\end{thm}
{\noindent{\bf Proof.}} In order to compute $v_{a\cdot\na f}$ we note that
\[
\Delta (a\cdot\na v_f) = \Delta a\cdot \na v_f + \na a\cdot\na\na v_f,
\]
and because $v_f\in e^{zl}H^{2}(Q)$ and $a\in B(\Omega)$ we have that
\[
\|\Delta (a\cdot\na v_f)\|_{L^2(e^{zl}dzdx)}\le C\|a\|_{B(\Omega)}\|v_f\|_{e^{zl}H^2(Q)}.
\]
Then solving
\[
\left \{
\ba
\Delta u = \Delta (a\cdot\na v_f) \quad {\mbox{in}}\; Q,\\
u = 0 \quad {\mbox{on}} \; \pa Q,
\ea
\right.
\]
we obtain $u\in H^{2}(Q)$ (by Lemma \ref{z}) and therefore $\pa_z u\in H^1_{0,L}(Q)$. Consequently  $-\pa_z u_{\left|\right. z=0}\in V_0(\Omega)$. Because $v_f$ vanishes on the boundary and $a\cdot\na$ is tangent to the boundary, it follows that $a\cdot\na v_f\in H_{0,L}^{1}(Q)$ (vanishes on the lateral boundary of $Q$ and is in $H^{1}(Q)$) and therefore
\[
v_{a\cdot\na f} = a\cdot\na v_f - u.
\]
Consequently
\[
[a\cdot\na, \l]f = -\pa_z u_{\left |\right. z=0}.
\]
The estimate (\ref{comtwo}) follows from the elliptic estimates and restriction estimates on $u$, as above:
\[
\|\pa_zu_{\left |\right. z=0}\|_{V_0}\le C \|\pa_z u\|_{H^1(Q)}
\le C\|a\|_{B(\Omega)}\|v_f\|_{e^{zl}H^2(Q)}\le C \|a\|_{B(\Omega)}\|f\|_{\fr{3}{2},D}
\]
\section{Linear transport and nonlocal diffusion}
We study the equation
\be
\pa_t\theta + u\cdot\na\theta + \l \theta = 0
\la{teq}
\ee
with initial data 
\be
\theta(x,0) = \theta_0
\la{it}
\ee
in the bounded open domain $\Omega\subset\Rr^d$ with smooth boundary. We assume that $u=u(x,t)$ is divergence-free
\be
\na\cdot u = 0,
\la{divuz}
\ee
that $u$ is smooth 
\be
u\in L^{2}(0,T; B(\Omega)^d),
\la{smu}
\ee
and that $u$ is parallel to the boundary 
\be
u_{\left |\right. \pa \Omega}\cdot n = 0.
\la{upa}
\ee
We consider zero boundary conditions for $\theta$. Strictly speaking, because this is a first order equation, it is better to think of these as a constraint on the evolution equation. We satrt with initial data $\theta_0$ which vanish on the boundary, and maintain this property in time. The transport evolution
\[
\pa_t\theta + u\cdot\na\theta = 0
\]
and, separately, the nonlocal diffusion
\[
\pa_t\theta + \l \theta = 0
\]
keep the constraint of $\theta_{\left |\right. \pa \Omega}= 0$. Because the
operators $u\cdot\na $ and $\l$ have the same differential order, neither dominates the other, and the linear evolution needs to be treated carefully. We start by considering Galerkin approximations. Let
\be
P_mf = \sum_{j=1}^m f_jw_j, \quad {\mbox{for}}\; f=\sum_{j=1}^{\infty}f_jw_j,
\la{pm}
\ee
and let
\be
\theta_m (x,t) = \sum_{j=1}^m\theta_j^{(m)}(t)w_j(x)
\la{thetam}
\ee
obey
\be
\pa_t\theta_m + P_m\left(u\cdot\na\theta_m\right) + \l\theta_m =0
\la{tetmeq}
\ee
with initial data
\be
\theta_m(x,0) = (P_m\theta_0)(x).
\la{thetmz}
\ee
These are ODEs for the coefficients $\theta_j^{(m)}(t)$, written conveniently. We prove bounds that are independent of $m$ and pass to the limit. Note that by construction
\[
\theta_m\in {\mathcal{D}}\left(\l^r\right), \quad \forall r\ge 0.
\]
We start with 
\be
\fr{1}{2}\fr{d}{dt}\|\theta_m\|^2_{L^2(\Omega)} + \|\theta_m\|^2_{V_0} = 0
\la{l2}
\ee
which implies
\be
\sup_{0\le t\le T}\fr{1}{2}\|\theta_m(\cdot, t)\|^2_{L^2(\Omega)} + \int_0^T\|\theta_m\|^2_{V_0}dt \le \fr{1}{2}\|\theta_0\|^2_{L^2(\Omega)}.
\la{l2b}
\ee
This follows because of the divergence-free condition and the fact that $u_{\left |\right. \pa\Omega}$ is parallel to the boundary. Next, we apply $\l$ to \eqref{tetmeq}. For convenience, we denote
\be
[\l, u\cdot\na]f = \Gamma f
\la{commu}
\ee
because $u$ is fixed throughout this section. Because $P_m$ and $\l$ commute, we have thus
\be
\pa_t\l\theta_m + P_m\left(u\cdot\na \l\theta_m + \G \theta_m\right) +\l^2\theta_m= 0.
\la{ltetmeq}
\ee
Now, we multiply \eqref{ltetmeq} by $\l^3\theta_m$ and integrate.  Note that
\[
\int_{\Omega} P_m\left(u\cdot\na \l\theta_m + \G \theta_m\right) \l^3\theta_m dx=
\int_{\Omega} \left(u\cdot\na \l\theta_m + \G \theta_m\right) \l^3\theta_m dx 
\]
because $P_m\theta_m = \theta_m$ and $P_m$ is selfadjoint. We bound the term
\[
\left |\int_{\Omega} \G\theta_m\l^3\theta_mdx\right |\le \|\G\theta_m\|_{V_0}\|\l^{2.5}\theta_m\|_{L^2(\Omega)}
\]
and use Theorem \ref{comthm} (\ref{comtwo}) to deduce
\[
\left |\int_{\Omega} \G\theta_m\l^3\theta_mdx\right |\le C\|u\|_{B(\Omega)}\|\l\theta_m\|_{V_0}\|\l^{2.5}\theta_m\|_{L^2(\Omega)}.
\]
We compute
\[
\ba
\int_{\Omega}(u\cdot\na\l\theta_m)\l^3\theta_m dx = \int_{\Omega}\l^2(u\cdot\na\l\theta_m)\l\theta_m \\
=\int_{\Omega}\left [(-\Delta u)\cdot\na \l\theta_m - 2\na u\cdot\na\na\l\theta_m\right]\l\theta_mdx  + \int_{\Omega}(u\cdot\na\l^3\theta_m)\l\theta_m dx\\
= \int_{\Omega}\left [(-\Delta u)\cdot\na \l\theta_m - 2\na u\cdot\na\na\l\theta_m\right]\l\theta_mdx - \int_{\Omega}\l^3\theta_m(u\cdot\na \l\theta_m)dx\\
=\int_{\Omega}\left [((-\Delta u)\cdot\na \l\theta_m)\l\theta_m + 2\na u\na\l\theta_m\na\l\theta_m\right]dx-\int_{\Omega}(u\cdot\na\l\theta_m)\l^3\theta_m dx.
\ea
\]
In the first integration by parts we used the fact that $\l^3\theta_m$ is a finite linear combination of eigenfunctions which vanish at the boundary. Then we use the fact that $\l^2 = -\Delta$ is local. In the last equality we integrated by parts using the fact that $\l \theta_m$ is a finite linear combination of eigenfunctions which vanish at the boundary and the fact that $u$ is divergence-free.
It follows that
\[
\int_{\Omega}(u\cdot\na\l\theta_m)\l^3\theta_m dx = \fr{1}{2}\int_{\Omega}\left [((-\Delta u)\cdot\na \l\theta_m)\l\theta_m + 2\na u\na\l\theta_m\na\l\theta_m\right]dx
\]
and consequently
\[
\left| \int_{\Omega}(u\cdot\na\l\theta_m)\l^3\theta_m dx\right| \le C\|u\|_{B(\Omega)}\|\l^2\theta_m\|_{L^2(\Omega)}^2
\]
We obtain thus
\be
\sup_{0\le t\le T}\|\l^2\theta_m(\cdot, t)\|^2_{L^2(\Omega)} + \int_0^T\|\l^2\theta_m\|^2_{V_0}dt \le C\|\l^2\theta_0\|^2_{L^2(\Omega)}e^{C\int_0^T\|u\|_{B(\Omega)}^2dt}.
\la{liv2b}
\ee
Passing to the limit $m\to\infty$ is done using the Aubin-Lions Lemma ({\cite{lions}}). We obtain
\beg{thm}{\la{linthm}} Let $u\in L^2(0,T; B(\Omega)^d)$ be a vector field parallel to the boundary. Then the equation (\ref{teq}) with initial data $\theta_0\in H_0^1(\Omega)\cap H^2(\Omega)$ has unique solutions belonging to
\[
\theta\in L^{\infty}(0,T; H^2(\Omega)\cap H_0^1(\Omega))\cap L^2(0,T; H^{2.5}(\Omega)).
\]
If the initial data $\theta_0\in L^p(\Omega)$, $1\le p\le \infty$, then
\be
\sup_{0\le t \le T}\|\theta(\cdot, t)\|_{L^p(\Omega)} \le \|\theta_0\|_{L^p(\Omega)}
\la{lp}
\ee
holds.
\end{thm}
The estimate (\ref{lp}) holds because, by use of Proposition \ref{cordoba} for the diffusive part and integration by parts for the transport part, we have for solutions of (\ref{teq})
\[
\fr{d}{dt}\|\theta\|^p_{L^p(\Omega)}\le 0,
\]
$1\le p<\infty$. The $L^{\infty}$ bound follows by taking the limit $p\to\infty$ in (\ref{lp}).
\section{SQG}
We consider now the equation
\be
\pa_t \theta + u\cdot\na\theta + \l\theta = 0
\la{sqg}
\ee
with 
\be
u = R_D^{\perp}\theta
\la{u}
\ee
and 
\be
R_D = \na\l^{-1}
\ee
in a bounded open domain in $\Omega\subset \Rr^2$ with smooth boundary.  Local existence of smooth solutions is possible to prove using methods similar to those developed above for linear drift-diffusion equations.
We will consider  weak solutions (solutions which satisfy the equations in the sense of distributions).

\beg{thm}{\la{weakgl}} Let $\theta_0\in L^2(\Omega)$ and let $T>0$. There exists a weak solution of (\ref{sqg})
\[
\theta\in L^{\infty}(0,T; L^2(\Omega))\cap L^2(0,T; V_0(\Omega))
\]
satisfying $lim_{t\to 0}\theta(t) = \theta_0$ weakly in $L^2(\Omega)$.
\end{thm}
\noindent{\bf Proof.} We consider Galerkin approximations, $\theta_m$ 
\[
\theta_m (x,t) = \sum_{j=1}^m \theta_j(t)w_j(x)
\]
obeying the ODEs (written conveniently as PDEs):
\[
\pa_t\theta_m + P_m\left[R_D^{\perp}(\theta_m)\cdot\na \theta_m\right] + \l \theta_m = 0
\]
with initial datum 
\[
\theta_m(0) = P_m(\theta_0).
\]
We observe that, multiplying by $\theta_m$ and integrating we have
\[
\fr{1}{2}\fr{d}{dt}\|\theta_m\|^2 + \|\theta_m\|^2_{\fr{1}{2}, D} = 0
\]
which implies that the sequence $\theta_m$ is bounded in
\[
\theta_m\in L^{\infty}(0,T; L^2(\Omega))\cap L^2(0,T; V_0(\Omega))
\]
It is known (\cite{cabre}) that $V_0(\Omega)\subset L^4(\Omega)$ with continuous inclusion. It is also known (\cite{jerisonkenig}) that 
\[
R_D: L^4(\Omega)\to L^4(\Omega)
\]
are bounded linear operators. It is then easy to see that 
$\pa_t \theta_m$ are bounded in $L^2(0,T; H^{-1}(\Omega))$. Applying the Aubin-Lions lemma, we obtain a subsequence, renamed $\theta_m$ converging strongly in $L^2(0,T; L^2(\Omega))$ and weakly in $L^2(0,T; V_0(\Omega))$ and in  $L^2(0,T; L^4(\Omega))$. The limit solves the equation (\ref{sqg}) weakly. Indeed, this follows after integration by parts because the product  $(R_D^{\perp}\theta_m)\theta_m$ is weakly convergent in $L^2(0,T; L^2(\Omega))$  by weak-times-strong weak continuity. The weak continuity in time at $t=0$ follows by integrating 
\[
(\theta_m(t), \phi) - (\theta_m(0), \phi) = \int_0^t\fr{d}{ds}\theta_m(s)ds
\]
and use of the equation and uniform bounds. We omit further details.\\

{\bf{Acknowledgment.}} The work of PC was partially supported by  NSF grant DMS-1209394

\end{document}